\documentclass[twoside,english]{paper}
\usepackage[T1]{fontenc}
\usepackage[latin9]{inputenc}
\usepackage{geometry}
\usepackage{babel}
\usepackage{mathrsfs}
\usepackage{amsmath}
\usepackage{amsthm}
\usepackage{amssymb}
\usepackage{booktabs}
\usepackage{enumitem}
\usepackage{longtable}
\usepackage{float}
\usepackage[pdfusetitle,
 bookmarks=true,bookmarksnumbered=false,bookmarksopen=false,
 breaklinks=false,pdfborder={0 0 1},backref=false,colorlinks=false]
 {hyperref}

\makeatletter
\let\paperclassexample\example

\let\example\relax
\AtBeginDocument{%
  \@ifundefined{example}{\let\example\paperclassexample}{}
}
\numberwithin{equation}{section}
\numberwithin{figure}{section}
\newcommand{\lyxaddress}[1]{
	\par {\raggedright #1
	\vspace{1.4em}
	\noindent\par}
}

\makeatother

\newtheorem{theorem}{Theorem}[section]
\newtheorem{proposition}[theorem]{Proposition}
\newtheorem{lemma}[theorem]{Lemma}
\newtheorem{corollary}[theorem]{Corollary}
\theoremstyle{definition}
\newtheorem{definition}[theorem]{Definition}
\newtheorem{example}[theorem]{Example}
\theoremstyle{remark}
\newtheorem{remark}[theorem]{Remark}

\newcommand{\OO}{\mathbb{O}}          
\newcommand{\MO}{\mathcal{O}}         
\newcommand{\RR}{\mathbb{R}}

\newcommand{\ZZ}{\mathbb{Z}}
\newcommand{\FF}{\mathbb{F}}

\newcommand{\LL}{\Lambda}             

\newcommand{\Nrm}{N}                  

\DeclareMathOperator{\Co}{Co}
\DeclareMathOperator{\Aut}{Aut}

\newcommand{\scriptref}[1]{\cite{repo}~%
  \href{https://github.com/koeplinger/leech\_alg/blob/main/python\_project/src/#1}%
  {\nolinkurl{python_project/src/#1}}}
\newcommand{\gapref}[1]{\cite{repo}~%
  \href{https://github.com/koeplinger/leech\_alg/blob/main/gap\_project/#1}%
  {\nolinkurl{gap_project/#1}}}

\begin{document}
\title{An order on the Leech lattice from a \texorpdfstring{$\ZZ_3$}{Z3}-symmetric
triple-octonion product}
\subtitle{1 November 2026 (v8)}
\author{Jens K\"oplinger}
\maketitle

\lyxaddress{105 E Avondale Dr, Greensboro, NC 27403, USA \\
jenskoeplinger@gmail.com}

\begin{abstract}
We show that the Leech lattice~$\LL$ admits the structure of an order
in a 24-dimensional real algebra: the lattice is closed under a
bilinear product on~$\RR^{24}$, assembled from nine
octonion-multiplication blocks under cyclic $\ZZ_3$ permutation of
three octonion factors.  In Robert A.~Wilson's sublattice description
of~$\LL$, closure turns on a transposition~$\sigma$ of two imaginary
basis units: $\sigma$ leaves the $E_8$ lattice~$L$ invariant and is
an octonion-algebra isomorphism, but moves Wilson's two sublattices
of~$L$.  The proof reduces to four lemmas, three of which
follow from the closure of one maximal order of integral octonions.  More
generally, for every linear map~$g$ with $g(L) = L$, the triple product
built from the $g$-twisted octonion product closes on~$\LL$ exactly when
Wilson's two sublattices are closed under that product.  We
determine the idempotents, the square-zero elements, and the automorphism
group of the algebra, as well as the automorphism group of the order and
the index of the span of its products.  We discuss the construction in the context of related work.
A historical appendix (1923--1946) traces
$\sigma$ to Richard H.~Bruck's repair of Johannes Kirmse's integral
octonions.  The
research record of this human--AI collaboration is public.
\end{abstract}



\section{Introduction}\label{sec:intro}

The Leech lattice is the unique even unimodular lattice of rank~24
with no vectors of squared length~2 \cite{ConwaySloane1999}.  Its
minimal shell consists of 196{,}560 vectors of squared length~4, and its
automorphism group is the Conway group~$\Co_0$, a double cover of the
sporadic simple group~$\Co_1$.  Throughout this paper, $\LL$ denotes
the $\sqrt{2}$-scaled copy of this lattice, at the scale of Wilson's
construction (Section~\ref{subsec:wilson-leech}); its minimal shell
$\operatorname{Min}(\LL)$ consists of vectors of squared length~8.

In rank~8, the $E_8$ lattice~$L$ is closed under an octonion
multiplication.  For a fixed multiplication, seven maximal orders of
integral octonions contain all eight basis units, each an $E_8$
lattice~\cite{Coxeter1946, ConwaySmith2003} (Appendix~\ref{app:history}).
Wilson~\cite[Sections~3--4]{Wilson2009} showed that $\LL$ embeds in~$L^3$ via
three sublattice conditions, raising the question: does
$\LL$ inherit an algebraic structure from the octonion product?

This question has been studied from several angles.  Geoffrey
M.~Dixon \cite{Dixon2010} obtained the integral octonions from an
octonion XY-product, represented~$\LL$ in~$\OO^3$, and asked whether
$\LL$ can be given its own algebraic structure~\cite[Section~7]{Dixon2010}.
John C.~Baez and Greg Egan \cite{Baez2014} explored
Jordan-algebraic structures containing the Leech lattice within the exceptional Jordan
algebra~$\mathfrak{h}_3(\OO)$, finding closure under a doubled
Jordan product.

Here we introduce a product on the direct sum
$\OO_1 \oplus \OO_2 \oplus \OO_3$ of three octonion blocks: it
multiplies blocks pairwise with the same octonion multiplication and
routes each product to one block by a $\ZZ_3$-symmetric rule.  Applied
to vectors of~$\LL$, the product closes once the octonion
multiplication is twisted by a transposition~$\sigma$ of two imaginary
basis units: we show that, with this twist, Wilson's sublattice
conditions yield the required closure inclusions.

\begin{theorem}\label{thm:main}
Let\/ $\star$ denote the $\ZZ_3$-symmetric triple-octonion product
on~$\RR^{24} = \OO_1 \oplus \OO_2 \oplus \OO_3$ of
Definition~\ref{def:triple}.  Then $u \star v \in \LL$ for all
$u, v \in \LL$ in Wilson's representation~\cite{Wilson2009}.
\end{theorem}

\begin{corollary}\label{cor:order}
$(\LL, +, \star)$ is an order\footnote{The ambient algebra has no
multiplicative identity (Remark~\ref{rem:no-identity}), so the term
``order'' is used here in the nonunital sense.} in the $\RR$-algebra
$(\RR^{24}, +, \star)$.
\end{corollary}

\begin{proof}
The product~$\star$ is bilinear (it is a sum of bilinear octonion
products).  By Theorem~\ref{thm:main}, $\star$ maps $\LL \times \LL$
into~$\LL$.
\end{proof}

\begin{remark}[Order-closure, not shell-closure]\label{rem:order-closure}
The closure $\LL \star \LL \subseteq \LL$ makes $(\LL, +, \star)$ an
order in the ring-theoretic sense: a subring that is a free
$\ZZ$-module of rank equal to the ambient algebra's dimension and spans
it.  The
product~$\star$ does not preserve the minimal shell
(Remark~\ref{rem:product-norms}).
Closure of the minimal shell, which spans~$\LL$, would be the stronger
property.
\end{remark}

The scope of this paper is deliberately bounded: we establish
existence of the order, with a proof that is elementary (four short
lemmas, three of which follow from the closure of one maximal
order).  Section~\ref{sec:prelim} fixes the octonion
multiplication, the $E_8$ lattice~$L$, the maximal order on which the
proof rests, and Wilson's description of~$\LL$.
Section~\ref{sec:construction} defines the product, and
Section~\ref{sec:proof} proves Theorem~\ref{thm:main}.
Section~\ref{sec:properties} determines algebraic properties of the
order: that standard identities fail on it; its idempotents; that it is
not maximal; its automorphism group; and its square-zero elements.  It
also characterizes the twists other than~$\sigma$ that close on~$\LL$.
Section~\ref{sec:related} places the construction among related work,
Section~\ref{sec:conclusion} concludes with the open questions, and
Section~\ref{sec:outlook} is an outlook.  Appendix~\ref{sec:appendix-tables}
reads the closure criterion modulo~$2$, Appendix~\ref{app:history} traces
the integral octonions from Leonard E.~Dickson to H.\,S.\,M.~Coxeter,
Appendix~\ref{sec:methodology} records the research methodology,
Appendix~\ref{app:samples} gives sampled frequencies of identities of
Section~\ref{sec:properties} and the norms of products of minimal
vectors, and Appendix~\ref{sec:census} tests the closure of other twists
and determines the index of the span of the products.

All results were developed through a systematic human--AI
collaboration, and the complete research record is public~\cite{repo}:
every prompt, computation, and verification, end to end
(Appendix~\ref{sec:methodology}).

\section{Preliminaries}\label{sec:prelim}

\subsection{The octonion algebra}

Let $\OO$ denote the real octonion algebra with basis
$\{e_0, e_1, \ldots, e_7\}$, where $e_0$ is the identity.
Its multiplication, the \emph{standard octonion product}, is
determined by the seven \emph{Fano triples}
\begin{equation}\label{eq:fano}
  (1,2,4),\; (2,3,5),\; (3,4,6),\; (4,5,7),\; (5,6,1),\;
  (6,7,2),\; (7,1,3),
\end{equation}
encoding the rule: for each triple $(a,b,c)$,
$e_a \cdot e_b = +e_c$, with cyclic permutations carrying a positive
sign and anti-cyclic permutations carrying a negative sign
\cite{Baez2002}.  Each imaginary basis unit $e_1, \ldots, e_7$ squares
to~$-e_0$.  Among
\emph{unital} real algebras, $\OO$ is the unique normed division
algebra of dimension~8 (Adolf Hurwitz); relaxing ``unital normed division
algebra'' to the absence of zero divisors (equivalently, bijectivity of left and
right multiplication by every nonzero element) admits further
examples, including the nonunital Okubo algebra.

The \emph{norm form} is $\Nrm(x) = x \bar{x} = \sum_{i=0}^7 x_i^2$
for $x = \sum_i x_i e_i$, where $\bar{x} = x_0 e_0 - \sum_{i=1}^7
x_i e_i$ is the conjugate.  Throughout, \emph{norm} means this
squared Euclidean length, and $\langle x, y\rangle$ is the Euclidean inner
product, so that $\langle x, x\rangle = \Nrm(x)$; on~$\RR^{24}$ both are
summed over the three blocks.  The composition property
$\Nrm(xy) = \Nrm(x)\Nrm(y)$ holds.  An octonion~$x$ is \emph{purely
imaginary} when $x_0 = 0$, that is, when $\bar x = -x$.

\subsection{The \texorpdfstring{$E_8$}{E8} lattice}

In Wilson's description~\cite{Wilson2009}, the $E_8$ lattice is
the set~$L \cong D_8^+$ of vectors in~$\RR^8$ whose coordinates (in the
basis $\{e_0, \ldots, e_7\}$) are either all integers with even
coordinate sum, or all half-integers with odd coordinate sum
\cite{ConwaySmith2003}.  It has
240~roots (vectors of norm~2).  Under the standard octonion
product, $L$~is closed: $L \cdot L \subseteq L$
(Lemma~\ref{lem:L-order}).  At this normalization $e_0 \notin L$ (the
lattice is even), so $L$~is a nonunital order.  As a lattice it is a
$\sqrt{2}$-scaled copy of a maximal order of integral octonions, the
unit-scale octavians, for which closure and maximality hold in the
classical unital sense \cite{Coxeter1946, ConwaySmith2003}.  For a
modern survey of the integral-octonion history,
see~\cite{Petersson2018}.

One maximal order of integral octonions is used throughout:
\begin{equation}\label{eq:MO}
  \MO := \tfrac12\,(e_1 + e_2)\,L .
\end{equation}
The lattice~$\MO$ contains the eight basis units and is closed under
the standard octonion product.  It is one of the seven maximal orders
of Section~\ref{sec:intro}, the ring $A_3$ of
Wilson~\cite[Section~2]{Wilson2009}, with Wilson's basis units $1$
and~$i_t$ ($0 \le t \le 6$) read as $e_0$ and~$e_{t+1}$.\footnote{\scriptref{verify\_structural\_proof.py}
\scriptref{verify\_sigma\_sublattice\_maximal\_order.py}.}  Since
$(e_1 + e_2)^2 = -2$ and $L = -L$, the left alternative law
$x(xy) = x^2 y$ gives $L = (e_1 + e_2)\,\MO$.

\subsection{Wilson's Leech lattice}\label{subsec:wilson-leech}

Define $s = \tfrac{1}{2}(-e_0 + e_1 + \cdots + e_7)$ and
$\bar{s} = \tfrac{1}{2}(-e_0 - e_1 - \cdots - e_7)$, both of
norm~$\Nrm(s) = \Nrm(\bar{s}) = 2$.  Then $s \in L$ (it is a
root in the half-integer branch, with an odd coordinate sum), while $\bar{s} \notin L$
(its coordinate sum is~$-4$, which falls in the excluded even
branch of the half-integer case).  Let
\begin{equation}\label{eq:disp1}
  L\bar{s} = \{ \ell \cdot \bar{s} : \ell \in L \}, \qquad
  Ls = \{ \ell \cdot s : \ell \in L \}.
\end{equation}
Both $L\bar{s}$ and~$Ls$ are sublattices of~$L$,\footnote{$L\bar
s \subseteq L$ even though $\bar s \notin L$: each $\ell \cdot
\bar s$ with $\ell$ a root of~$L$ lies in~$L$ by direct
computation, and the inclusion extends by $\ZZ$-linearity in the
left factor.  $L\bar s$ is therefore a $\ZZ$-subgroup of~$L$.
For~$Ls$, the same follows from $s \in L$ and
$L \cdot L \subseteq L$ (Lemma~\ref{lem:L-order}).}
with
$[L : L\bar{s}] = [L : Ls] = 16$,
$L\bar{s} + Ls = L$, and $L\bar{s} \cap Ls = 2L$.

\begin{definition}[Wilson]\label{def:wilson-leech}
The \emph{Leech lattice} is
\begin{equation}\label{eq:disp2}
  \LL = \bigl\{ (x,y,z) \in L^3 :
    \text{(1) } x,y,z \in L;\;
    \text{(2) } x{+}y,\, x{+}z,\, y{+}z \in L\bar{s};\;
    \text{(3) } x{+}y{+}z \in Ls
  \bigr\}.
\end{equation}
\end{definition}

The 196{,}560 minimal vectors of~$\LL$ decompose into three families:

\smallskip
\begin{tabular}{@{}lll@{}}
\toprule
Family & Formula & Count \\
\midrule
Type~1 & $(2\lambda, 0, 0)$ and cyclic permutations & $3 \times 240 = 720$ \\
Type~2 & $(\lambda\bar{s},\; (\lambda\bar{s})j,\; 0)$ and cyclic
  permutations & $3 \times 240 \times 16 = 11{,}520$ \\
Type~3 & $((\lambda s)j,\; \lambda k,\; (\lambda j)k)$ and cyclic
  permutations & $3 \times 240 \times 16^2 = 184{,}320$ \\
\bottomrule
\end{tabular}
\smallskip

\noindent
Here $\lambda$ ranges over the 240~roots of~$L$, and
$j, k \in J = \{\pm e_t : 0 \le t \le 7\}$ (16~unit octonions).

\section{The construction}\label{sec:construction}

\begin{definition}[Transposition]\label{def:twist}
Let $\sigma = (i\; j)$ be a transposition on $\{1, \ldots, 7\}$,
extended to the linear involution $\sigma \colon \RR^8 \to \RR^8$
that fixes $e_0$ and swaps $e_i \leftrightarrow e_j$.
\end{definition}

\begin{definition}[$\sigma$-twisted octonion product]\label{def:twisted-product}
The \emph{$\sigma$-twisted octonion product}~$\cdot_\sigma$ on~$\OO$
is defined from the standard octonion product~\eqref{eq:fano} by
\begin{equation}\label{eq:twist-formula}
  x \cdot_\sigma y \;:=\; \sigma\bigl(\sigma(x) \cdot \sigma(y)\bigr).
\end{equation}
\end{definition}

This is another octonion product on $\OO$: by definition,
$\sigma \colon (\OO, \cdot) \to (\OO, \cdot_\sigma)$ is an algebra
isomorphism, $\sigma(x \cdot y) = \sigma(x) \cdot_\sigma \sigma(y)$.

\begin{definition}[Triple product with $\ZZ_3$ routing]
\label{def:triple}
Decompose $\RR^{24} = \OO_1 \oplus \OO_2 \oplus \OO_3$ into three
copies of~$\RR^8$ (blocks indexed $\alpha = 1, 2, 3$).  For
$u = (x, y, z)$ and $v = (x', y', z')$ with $x, x' \in \OO_1$, etc.,
define the \emph{$\ZZ_3$-symmetric triple-octonion product}, or
\emph{triple product} for short,
\begin{equation}\label{eq:triple}
  u \star v = (P_1,\, P_2,\, P_3),
\end{equation}
where each block $P_\alpha$ is the sum of three $\sigma$-twisted
octonion products:
\begin{align}
  P_1 &= x \cdot_\sigma x' + z \cdot_\sigma y' + y \cdot_\sigma z',
  \label{eq:P1} \\
  P_2 &= y \cdot_\sigma y' + x \cdot_\sigma z' + z \cdot_\sigma x',
  \label{eq:P2} \\
  P_3 &= z \cdot_\sigma z' + y \cdot_\sigma x' + x \cdot_\sigma y'.
  \label{eq:P3}
\end{align}
The routing rule is: a product of blocks $\alpha$ and~$\beta$ lands in
block~$\gamma$, where $\{\alpha, \beta, \gamma\} = \{1,2,3\}$ for
$\alpha \ne \beta$, and in block~$\alpha$ when $\alpha = \beta$.
All nine block-pair products use the same multiplication~$\cdot_\sigma$.
\end{definition}

\section{Proof of closure}\label{sec:proof}

We prove Theorem~\ref{thm:main} by showing that the $\sigma$-twisted
octonion product preserves each of Wilson's three conditions
(Definition~\ref{def:wilson-leech}).  Throughout, and without loss of
generality, $\sigma = (1\;2)$.\footnote{The 21 transpositions of
imaginary basis units form a single orbit under conjugation by the
order-21 group of coordinate permutations preserving the oriented Fano
triples~\eqref{eq:fano}.  Each such permutation is an octonion
automorphism that fixes $s$ and~$\bar s$ and preserves~$L$, so
conjugating by it carries the proof below for $(1\;2)$ to any other
transposition.}  The map~$\sigma$ fixes~$L$ setwise
(Lemma~\ref{lem:sigma-L}) and carries the standard octonion product to
an isomorphic one (Definition~\ref{def:twisted-product}).  Its
load-bearing effect is that it \emph{moves} the sublattices $L\bar s$
and~$Ls$ onto left multiples of the maximal order~$\MO$
of~\eqref{eq:MO}: $\sigma(L\bar s) = 2\MO$ and $\sigma(Ls) = p\MO$ for a
purely imaginary $p \in \MO$.  Lemmas~\ref{lem:L-sigmaLsbar}
and~\ref{lem:sigmaLs-closed}, the closure inclusions for the moved
sublattices, then follow from the closure of~$\MO$.

\subsection{The four lemmas}

\begin{lemma}\label{lem:sigma-L}
$\sigma(L) = L$.
\end{lemma}

\begin{proof}
The transposition~$\sigma$ permutes two coordinates, so it preserves
the two conditions that define~$L$ (Section~\ref{sec:prelim}):
all coordinates integers or all half-integers, and the parity of the
coordinate sum.  Hence
$\sigma(L) \subseteq L$.  Since $\sigma$ is an involution,
$\sigma(L) = L$.
\end{proof}

The remaining three lemmas rest on the maximal order~$\MO$
of~\eqref{eq:MO} and on three facts.  First, two lattice identities:
with $p := e_2 - e_4 + e_5 - e_6$, a purely imaginary element of~$\MO$
of norm~$4$,
\begin{equation}\label{eq:multiples}
  \sigma(L\bar s) = 2\MO, \qquad \sigma(Ls) = p\MO .
\end{equation}
The inclusions $2\MO \subseteq \sigma(L\bar s)$ and
$p\MO \subseteq \sigma(Ls)$ are checked on a $\ZZ$-basis $\{b_k\}$
of~$\MO$: the eight vectors $2b_k$ lie in $\sigma(L\bar s)$, and the
eight products $p\,b_k$ lie in $\sigma(Ls)$.  Both sides of each
identity have covolume~$16$.  Left multiplication by~$x$ scales inner
products by $\Nrm(x)$ and hence covolumes by $\Nrm(x)^4$.  Since
$L = (e_1 + e_2)\,\MO$ has covolume~$1$ ($E_8$ is unimodular), $\MO$
has covolume~$\tfrac1{16}$, so $2\MO$ and $p\MO$ have covolume~$16$.
So do $\sigma(L\bar s)$ and~$\sigma(Ls)$, by
$[L : L\bar s] = [L : Ls] = 16$ (Section~\ref{subsec:wilson-leech}).
A sublattice of the same covolume is the whole
lattice.\footnote{\scriptref{verify\_structural\_proof.py}.}

Second, an identity in the octonions: for purely imaginary~$x$ and all
$a, b \in \OO$,
\begin{equation}\label{eq:moufang-imaginary}
  (xa)(xb) = x\bigl((ax)b\bigr).
\end{equation}
Indeed, both sides vanish for $x = 0$.  For $x \ne 0$ we have
$x^2 = -x\bar x = -\Nrm(x)$, so $a = yx$ with $y := -ax/\Nrm(x)$ by the
right alternative law $(ax)x = ax^2$, and $xa = xyx$ by the flexible
law $x(yx) = (xy)x$.  The left Moufang identity $(xyx)z = x(y(xz))$
\cite{ConwaySmith2003} with $z = xb$ and the left alternative law
$x(xb) = x^2 b$ give
$(xa)(xb) = x\bigl(y(x^2 b)\bigr) = -\Nrm(x)\,x(yb) = x\bigl((ax)b\bigr)$.

Third, its consequence for lattices: if $M \subseteq \OO$ is closed
under the standard octonion product and $x \in M$ is purely imaginary,
then
by~\eqref{eq:moufang-imaginary}
\begin{equation}\label{eq:multiple-closed}
  (xM)(xM) \subseteq x\bigl((Mx)M\bigr) \subseteq xM .
\end{equation}

\begin{lemma}\label{lem:L-order}
$L \cdot L \subseteq L$.
\end{lemma}

\begin{proof}
$L = (e_1 + e_2)\,\MO$ by~\eqref{eq:MO}, with $e_1 + e_2$ a purely
imaginary element of~$\MO$; apply~\eqref{eq:multiple-closed} with
$M = \MO$.  This is the norm-2-scale form of the classical fact that
the maximal orders of integral octonions are closed under
multiplication \cite{Coxeter1946, ConwaySmith2003}.
\end{proof}

\begin{lemma}\label{lem:L-sigmaLsbar}
$L \cdot \sigma(L\bar{s}) \subseteq \sigma(L\bar{s})$.
\end{lemma}

\begin{proof}
$L = (e_1 + e_2)\,\MO \subseteq \MO \cdot \MO \subseteq \MO$, so
by~\eqref{eq:multiples}
$L \cdot \sigma(L\bar s) = L \cdot 2\MO \subseteq 2\,(\MO \cdot \MO)
\subseteq 2\MO = \sigma(L\bar s)$.
\end{proof}

\begin{lemma}\label{lem:sigmaLs-closed}
$\sigma(Ls) \cdot \sigma(Ls) \subseteq \sigma(Ls)$.
\end{lemma}

\begin{proof}
$\sigma(Ls) = p\MO$ by~\eqref{eq:multiples}, with $p$ a purely
imaginary element of~$\MO$; apply~\eqref{eq:multiple-closed} with
$M = \MO$.
\end{proof}

\begin{remark}\label{rem:Ls-not-closed}
The sublattice~$Ls$ is \emph{not} closed under the standard octonion
product: $u_1 = (-1, 0, 0, 0, 1, 0, 1, 1)$ and
$u_2 = (-1, 1, 0, 0, 0, 1, 0, 1)$ lie in~$Ls$;
$u_1 \cdot u_2 = (0, -2, 2, 1, -2, -1, -1, -1)$ lies in
$L\bar{s} \setminus Ls$.  Moreover $\sigma(Ls) \ne Ls$: with
$\sigma = (1\;2)$, the
element $v = (-1, 0, 1, 0, 0, 1, 0, 1)$ lies in~$\sigma(Ls)
\setminus Ls$.  Likewise $\sigma(L\bar s) \ne
L\bar s$: $w = (0, 0, -1, 0, -1, 0, -1, -1) \in \sigma(L\bar s)
\setminus L\bar s$.  The closure of~$\sigma(Ls)$ is not inherited
from an ideal property: $\sigma(Ls) = p\MO$ is not an ideal of~$L$ on
either side, since neither $L \cdot \sigma(Ls)$ nor
$\sigma(Ls) \cdot L$ lies
in~$\sigma(Ls)$.\footnote{\scriptref{verify\_sigma\_Ls\_ideal\_exclusion.py}.}
By contrast, $\sigma(L\bar s) = 2\MO$ is a two-sided ideal of~$L$, of
which Lemma~\ref{lem:L-sigmaLsbar} is the left-handed half.
\end{remark}

\subsection{Blockwise form of the triple product}

Since each $\sigma$-twisted product has the form
$a \cdot_\sigma b = \sigma(\sigma(a) \cdot \sigma(b))$ and $\sigma$
is linear, each block $P_\alpha$ is $\sigma$ applied to a sum of
standard octonion products of $\sigma$-images of the inputs.
Write $X = \sigma(x)$, $Y = \sigma(y)$, $Z = \sigma(z)$,
$X' = \sigma(x')$, $Y' = \sigma(y')$, $Z' = \sigma(z')$:
\begin{align}
  P_1 &= \sigma\bigl(X \cdot X' + Z \cdot Y' + Y \cdot Z'\bigr),
  \label{eq:P1-std} \\
  P_2 &= \sigma\bigl(Y \cdot Y' + X \cdot Z' + Z \cdot X'\bigr),
  \label{eq:P2-std} \\
  P_3 &= \sigma\bigl(Z \cdot Z' + Y \cdot X' + X \cdot Y'\bigr).
  \label{eq:P3-std}
\end{align}
Taking sums:
\begin{itemize}
\item \textbf{Block sum.}
  $P_1 + P_2 + P_3 = \sigma\bigl((X+Y+Z) \cdot (X'+Y'+Z')\bigr)$.
\item \textbf{Pairwise sums.}  For example,
  $P_1 + P_2 = \sigma\bigl(X \cdot (X'+Z') + Y \cdot (Y'+Z') +
  Z \cdot (X'+Y')\bigr)$, and similarly for the other two pairs.
\end{itemize}

\subsection{Proof of Wilson's three conditions}

\subsubsection{Condition~(1): each block lies in~\texorpdfstring{$L$}{L}}

\begin{proposition}\label{prop:cond1}
For all $u, v \in \LL$, each block $P_\alpha$ of $u \star v$ lies
in~$L$.
\end{proposition}

\begin{proof}
Consider $P_1 = \sigma(X \cdot X' + Z \cdot Y' + Y \cdot Z')$.
Since $u \in \LL$, condition~(1) gives $x \in L$, so
$X = \sigma(x) \in L$ by Lemma~\ref{lem:sigma-L}.  Similarly
$X', Y, Y', Z, Z' \in L$.
By Lemma~\ref{lem:L-order}, each product $X \cdot X'$, $Z \cdot Y'$,
$Y \cdot Z'$ lies in~$L$.  Their sum lies in~$L$ (closed under
addition).  Applying $\sigma$ and using Lemma~\ref{lem:sigma-L},
$P_1 \in L$.  The same argument applies to
$P_2$ and~$P_3$.
\end{proof}

\subsubsection{Condition~(2): pairwise sums lie
in~\texorpdfstring{$L\bar{s}$}{L s-bar}}

\begin{proposition}\label{prop:cond2}
For all $u, v \in \LL$, the pairwise sums $P_\alpha + P_\beta$ lie
in~$L\bar{s}$.
\end{proposition}

\begin{proof}
We show $P_1 + P_2 \in L\bar{s}$; the other two cases follow by the
$\ZZ_3$ symmetry of the routing.

From \eqref{eq:P1-std}--\eqref{eq:P2-std}:
\begin{equation}\label{eq:disp3}
  P_1 + P_2 = \sigma\bigl(X \cdot (X' + Z') + Y \cdot (Y' + Z')
    + Z \cdot (X' + Y')\bigr).
\end{equation}

Since $v \in \LL$, condition~(2) gives $x' + z' \in L\bar{s}$, so
$X' + Z' = \sigma(x' + z') \in \sigma(L\bar{s})$.  Similarly
$Y' + Z' \in \sigma(L\bar{s})$ and $X' + Y' \in \sigma(L\bar{s})$.
Since $u \in \LL$, we have $X, Y, Z \in L$ (by Lemma~\ref{lem:sigma-L}).

By Lemma~\ref{lem:L-sigmaLsbar}, each product $X \cdot (X'+Z')$, $Y \cdot (Y'+Z')$,
$Z \cdot (X'+Y')$ lies in~$\sigma(L\bar{s})$.  Their sum lies
in~$\sigma(L\bar{s})$.  Applying the outer~$\sigma$:
\begin{equation}\label{eq:disp4}
  P_1 + P_2 \in \sigma(\sigma(L\bar{s})) = L\bar{s},
\end{equation}
since $\sigma^2 = \mathrm{id}$.
\end{proof}

\subsubsection{Condition~(3): the block sum lies
in~\texorpdfstring{$Ls$}{Ls}}

This is the crux of the proof.  With the standard octonion product in
place of~$\cdot_\sigma$, the block sum would require
$Ls \cdot Ls \subseteq Ls$, which does
not hold (Remark~\ref{rem:Ls-not-closed}).  The twist redirects the
computation through~$\sigma(Ls)$, where closure holds.

\begin{proposition}\label{prop:cond3}
For all $u, v \in \LL$, the block sum $P_1 + P_2 + P_3$ lies in~$Ls$.
\end{proposition}

\begin{proof}
From \eqref{eq:P1-std}--\eqref{eq:P3-std}:
\begin{equation}\label{eq:disp5}
  P_1 + P_2 + P_3 = \sigma\bigl((X+Y+Z) \cdot (X'+Y'+Z')\bigr).
\end{equation}

Since $u \in \LL$, condition~(3) gives $x + y + z \in Ls$, so
\begin{equation}\label{eq:disp6}
  X + Y + Z = \sigma(x+y+z) \in \sigma(Ls).
\end{equation}
Similarly $X' + Y' + Z' \in \sigma(Ls)$.

By Lemma~\ref{lem:sigmaLs-closed}, $\sigma(Ls) \cdot \sigma(Ls) \subseteq \sigma(Ls)$, so
\begin{equation}\label{eq:disp7}
  (X+Y+Z) \cdot (X'+Y'+Z') \in \sigma(Ls).
\end{equation}
Applying the outer~$\sigma$:
\begin{equation}\label{eq:disp8}
  P_1 + P_2 + P_3 \in \sigma(\sigma(Ls)) = Ls. \qedhere
\end{equation}
\end{proof}

This completes the proof of Theorem~\ref{thm:main}.

\section{Algebraic properties}\label{sec:properties}

Throughout this section, $\sigma = (1\;2)$, as in
Section~\ref{sec:proof}.  Appendix~\ref{sec:census} complements this
section: it tests other twists (Section~\ref{sec:criterion}) for
closure on~$\LL$ and determines how much of~$\LL$ the products of the
order reach.

\subsection{Identities}\label{sec:identities}

Several standard algebraic identities, among them associativity,
alternativity, flexibility, power-associativity, and commutativity,
fail on the ambient algebra $(\RR^{24}, +, \star)$ and hence, by
scaling, on the order $(\LL, +, \star)$.  The table below lists a
witness of each violation, which is easily checked by hand from
Definitions~\ref{def:twisted-product}
and~\ref{def:triple}.\footnote{\scriptref{verify\_identity\_witnesses.py}.}
Appendix~\ref{app:samples} estimates by sampling how often the
identities from left alternativity to symmetric composition hold on
minimal vectors.

\begin{center}
\begin{tabular}{@{}ll@{}}
\toprule
Identity & Witness of violation $(u, v\,[, w])$ \\
\midrule
Associativity $(u\star v)\star w = u\star(v\star w)$ &
  $(e_0, 0, 0)$, $(0, e_0, 0)$, $(0, 0, e_0)$ \\
Left alternativity $(u\star u)\star v = u\star(u\star v)$ &
  $(e_0, 0, 0)$, $(0, e_0, 0)$ \\
Right alternativity $(u\star v)\star v = u\star(v\star v)$ &
  $(0, e_0, 0)$, $(e_0, 0, 0)$ \\
Flexibility $(u\star v)\star u = u\star(v\star u)$ &
  $(e_3, e_4, 0)$, $(e_0, 0, 0)$ \\
Cube power-associativity $(u\star u)\star u = u\star(u\star u)$ &
  $(e_0, e_3, e_4)$ \\
Quartic power-associativity: the five parenthesizations &
  $(e_0, -e_0, 0)$ \\
\quad of $u\star u\star u\star u$ coincide & \\
Commutativity $u\star v = v\star u$ &
  $(e_3, 0, 0)$, $(e_4, 0, 0)$ \\
Norm multiplicativity $\Nrm(u\star v) = \Nrm(u)\Nrm(v)$ &
  $(e_0, e_0, e_0)$, $(e_0, -e_0, 0)$ \\
Symmetric composition $\langle u\star v, w\rangle =
  \langle u, v\star w\rangle$ &
  $(e_0, 0, 0)$, $(e_3, 0, 0)$, $(e_3, 0, 0)$ \\
Multiplicative identity &
  no identity element exists (Remark~\ref{rem:no-identity}) \\
\bottomrule
\end{tabular}
\end{center}

The failure of norm multiplicativity rules out $(\LL, +, \star)$
being a composition algebra for~$\Nrm$.  The failure of
symmetric composition rules out, independently, its being a symmetric
composition algebra for~$\Nrm$, such as the Okubo and para-Hurwitz
algebras (Alberto Elduque~\cite{Elduque2000_Triality,
Elduque2023_IsotropicNorm}; Susumu Okubo~\cite{Okubo1995_Book}).

\begin{remark}\label{rem:product-norms}
The product does not preserve the minimal shell: for the minimal
vector $u = (2\lambda, 0, 0)$ with $\lambda = e_3 + e_4$ (Type~1 in
Section~\ref{subsec:wilson-leech}), $u \star u = (4\lambda^2, 0, 0) =
(-8e_0, 0, 0)$, of norm~64.  Appendix~\ref{app:samples} records the
self-product norms of all minimal vectors and the pair-product norms on
sampled pairs.
\end{remark}

\begin{remark}\label{rem:no-identity}
The ambient algebra $(\RR^{24}, +, \star)$ has no multiplicative
identity.  A left identity~$(a,b,c)$ would act on~$(x',y',z')$ by the
block matrix
$\bigl(\begin{smallmatrix} a & c & b \\ c & b & a \\ b & a & c
\end{smallmatrix}\bigr)$, with entries acting by left multiplication
in $(\OO, \cdot_\sigma)$; each of $a$, $b$, $c$ occurs both on and off
the diagonal, so it would have to act both as the identity and as
zero.  In particular, $(e_0, e_0, e_0)$ is
not an identity: it maps $(x', y', z')$ to $(x'+y'+z',\, x'+y'+z',\,
x'+y'+z')$.  Consequently $(\LL, +, \star)$ is a nonunital order.
\end{remark}

The defining property of the order is \emph{closure}: the Leech
lattice is stable under the product.

\subsection{Idempotents}\label{sec:idempotents}

The ambient algebra $(\RR^{24}, +, \star)$ has exactly eight
idempotents, by symbolic
classification.\footnote{\scriptref{verify\_idempotent\_classification.py}.
The nonzero idempotents of~\eqref{eq:idempotents} are also obtained
in~\cite{Petersson2026}.}
Writing $\varepsilon_1 = (e_0, 0, 0)$, $\varepsilon_2 = (0, e_0, 0)$,
$\varepsilon_3 = (0, 0, e_0)$, and $\omega = \varepsilon_1 +
\varepsilon_2 + \varepsilon_3$, they are
\begin{equation}\label{eq:idempotents}
  0, \qquad \varepsilon_1, \ \varepsilon_2, \ \varepsilon_3, \qquad
  \tfrac13 \omega, \qquad
  \varepsilon_i - \tfrac13 \omega \ \ (i = 1, 2, 3).
\end{equation}
Note that every block of every one of them is a multiple of~$e_0$.

The order $(\LL, +, \star)$ has no nonzero
idempotent:\footnote{\scriptref{verify\_idempotent\_lattice\_rescaling.py}.}
each nonzero idempotent spans a rational ray that meets~$\LL$, and its
least positive multiple in~$\LL$ is
\begin{equation}\label{eq:idempotent-multiples}
  4\varepsilon_i \ (\Nrm = 16), \qquad
  4\omega \ (\Nrm = 48), \qquad
  6\bigl(\varepsilon_i - \tfrac13\omega\bigr) \ (\Nrm = 24).
\end{equation}
These multiples follow from Wilson's three conditions.  Since $L$,
$L\bar s$, and~$Ls$ meet~$\RR e_0$ in $2\ZZ e_0$, $2\ZZ e_0$,
and~$4\ZZ e_0$, a real multiple of an idempotent lies in~$\LL$ exactly
when its blocks lie in~$2\ZZ e_0$ (condition~(1), which then implies
condition~(2)) and its block sum lies in~$4\ZZ e_0$ (condition~(3)).
The intersections with $L\bar s$ and~$Ls$ follow from
$2e_0 \in \sigma(2\MO) = L\bar s$ (by~\eqref{eq:multiples}, as $\sigma$
fixes $e_0 \in \MO$) and from $L\bar s \cap Ls = 2L$, which contains
$4e_0$ but not~$2e_0$.
On the line spanned by a nonzero idempotent, the points of~$\LL$ are
the integer multiples $kw$ of the corresponding vector~$w$
in~\eqref{eq:idempotent-multiples}, of norm $16k^2$, $48k^2$, or $24k^2$;
in particular, no real multiple of an idempotent lies on the minimal
shell.

\subsection{Maximality}\label{sec:maximality}

The order $(\LL, +, \star)$ is not maximal: larger orders in the
ambient algebra contain it.  Two such orders, $L^3$ and
$\LL + \ZZ\varepsilon_1$, witness this, and neither contains the
other.

The first witness is~$L^3$.  The proof of
Proposition~\ref{prop:cond1} uses only Wilson's condition~(1), that
the blocks of both factors lie in~$L$, so it shows $L^3 \star L^3
\subseteq L^3$.  The lattice~$L^3$ is even and integral, and
$[L^3 : \LL] = 2^{12}$.

The second witness is $\LL + \ZZ\varepsilon_1$, the lattice generated
by~$\LL$ and the idempotent $\varepsilon_1 = (e_0, 0, 0)$.  On either
side, $\varepsilon_1$ acts as the swap of blocks~2 and~3: since $e_0$
is the identity of $(\OO, \cdot_\sigma)$, \eqref{eq:P1}--\eqref{eq:P3}
give $\varepsilon_1 \star v = v \star \varepsilon_1 = (x', z', y')$ for
$v = (x', y', z')$.  This swap preserves~$\LL$, because Wilson's three
conditions are symmetric in the blocks.  With $\varepsilon_1 \star
\varepsilon_1 = \varepsilon_1$, the lattice $\LL + \ZZ\varepsilon_1$ is
therefore $\star$-closed by bilinearity.
By~\eqref{eq:idempotent-multiples}, $4\varepsilon_1$ is the least
positive multiple of~$\varepsilon_1$ in~$\LL$, so
$[\LL + \ZZ\varepsilon_1 : \LL] = 4$.  Unlike~$L^3$, the second witness
is not integral: $\varepsilon_1$ has half-integer inner products with
some vectors of~$\LL$.

Neither witness contains the other: $\LL + \ZZ\varepsilon_1$ is not
contained in~$L^3$, because the first block~$e_0$ of~$\varepsilon_1$ is
not in~$L$, and $L^3$ is not contained in~$\LL + \ZZ\varepsilon_1$,
since $[L^3 : \LL] = 2^{12}$ exceeds $[\LL + \ZZ\varepsilon_1 : \LL] =
4$.

\subsection{The automorphism group}\label{sec:autgroup}

Restricted to the Leech lattice, the ambient algebra has the
automorphism group $\Aut(\LL, +, \star)$, the group of bijections
of~$\LL$ that preserve both $+$ and~$\star$.  This group is finite, of
order~$36$, with structure $C_6 \times S_3$.  The~$S_3$ permutes the three blocks of
Definition~\ref{def:triple}, and the~$C_6$ is generated by the
following signed permutation of imaginary basis units, an automorphism of
the algebra $(\OO, \cdot_\sigma)$ of
Definition~\ref{def:twisted-product}, applied in each block, with
$e_0$ fixed:
\begin{align}
  &e_1 \mapsto -e_5, \quad e_5 \mapsto e_6, \quad e_6 \mapsto e_1,
    \label{eq:c6-gen-1} \\
  &e_2 \mapsto -e_7, \quad e_7 \mapsto -e_3, \quad e_3 \mapsto e_2,
    \label{eq:c6-gen-2} \\
  &e_4 \mapsto -e_4. \label{eq:c6-gen-3}
\end{align}

Every automorphism of the order $(\LL, +, \star)$ is an isometry
of~$\LL$, so $\Aut(\LL, +, \star)$ is a subgroup of the Conway
group~$\Co_0$, the isometry group of~$\LL$.  This inclusion can already
be settled at the level of the ambient algebra, to which every
automorphism of the order extends linearly ($\LL$ spans~$\RR^{24}$, and
$\star$ is bilinear).  As a real algebra, $\RR^{24} = \OO^3$
splits as $D \oplus T$, where $D = \{(a,a,a) : a \in \OO\}$ (of
dimension~8) and $T = \{(x,y,z) \in \OO^3 : x+y+z = 0\}$ (of
dimension~16).  Both are two-sided ideals, and $D \star T = T \star D
= 0$.  Every automorphism of the ambient algebra preserves them and
acts on each by an octonion automorphism, composed with a permutation
of the blocks, so that the ambient algebra itself has the automorphism
group $\Aut(\RR^{24}, \star) = G_2 \times G_2 \times
S_3$, with $G_2$ the automorphism group of the
octonions.\footnote{\cite{Petersson2026}.}  An automorphism of
the \emph{real} octonions preserves the positive-definite octonion norm,
and $D \perp T$; hence every automorphism of the ambient algebra is
orthogonal, and
\begin{equation}\label{eq:disp9}
  \Aut(\LL, +, \star) \;\subsetneq\; \Co_0 .
\end{equation}
The inclusion is strict because $-I_{24} \in \Co_0$ does not
preserve~$\star$.

The derivation, the complete enumerations, and an
independent verification are recorded separately in a companion
note.\footnote{\cite{repo}~%
\href{https://github.com/koeplinger/leech\_alg/blob/main/paper/automorphism\_group\_2026-07-12.pdf}%
{\texttt{paper/automorphism\_group\_2026-07-12.pdf}};
\scriptref{verify\_aut\_lambda\_star.py}
\gapref{aut\_lambda\_star.g}.}

\subsection{Square-zero elements}\label{sec:squarezero}

In the ambient algebra, $u = (x, y, z)$ satisfies $u \star u = 0$ if
and only if $x$, $y$, $z$ are purely imaginary, of equal norm, and sum
to zero (again by symbolic classification).  Three vectors of equal norm summing to zero lie at $120^\circ$
to one another; with their negatives they are the six vertices of a
regular hexagon, and $x$, $y$, $z$ are three alternating vertices, a
\emph{hexagonal triple}.  The solution set is
a cone of dimension~12, and it lies in~$T$.\footnote{\scriptref{verify\_square\_zero\_classification.py}.}

The lattice does contain such elements.  A vector $(x, y, z) \in \LL$
squares to zero exactly when $x$, $y$, $z$ are purely imaginary vectors
of~$L\bar s$, have equal norms, and sum to zero (for $x + y + z = 0$,
Wilson's three conditions reduce to $x, y, z \in L\bar s$).  The 126 purely imaginary
vectors of norm~4 in~$L\bar s$ form an $E_7$ root system:
by~\eqref{eq:multiples}, $L\bar s = \sigma(2\MO)$ is an $E_8$ lattice
whose minimal vectors have norm~4 and include~$2e_0$, and these are its
minimal vectors orthogonal to~$2e_0$.  Each of these roots has exactly
32 roots at $120^\circ$; the third block is then
forced, so the ordered hexagonal triples number $126 \times 32 = 4{,}032$,
one square-zero vector of~$\LL$ each, all of norm~12.  The norm of
every square-zero vector of~$\LL$ lies in~$12\ZZ$, since
$\Nrm(u) = 3\Nrm(x)$ and every norm of~$L\bar s$ lies in~$4\ZZ$.  These
$4{,}032$ vectors are exactly the square-zero vectors of norm~12, and
norm~24 occurs as well: the sum of two of them that are blockwise
orthogonal is square-zero of norm~24.\footnote{\scriptref{verify\_square\_zero\_norm24\_witness.py}.}

\subsection{The closure criterion}\label{sec:criterion}

The construction twists the octonion product by one
transposition,~$\sigma$.  This subsection replaces~$\sigma$ by other
maps and characterizes those for which the resulting triple product
still closes on~$\LL$.

A \emph{twist} is an invertible linear map~$g$ of~$\RR^8$ with
$g(L) = L$.  It defines the twisted product $x \cdot_g y :=
g^{-1}\bigl(g(x) \cdot g(y)\bigr)$ on~$\OO$.  The triple
product~$\star_g$ is that of Definition~\ref{def:triple}, with
$\cdot_g$ in all nine block-pair products.  Since $\sigma^{-1} = \sigma$,
the twist $g = \sigma$ (Lemma~\ref{lem:sigma-L}) gives the product of
Definition~\ref{def:twisted-product}, and $\star_\sigma = \star$.

\begin{proposition}\label{prop:criterion}
For every twist~$g$, the product~$\star_g$ closes on~$\LL$ if and only
if $L\bar s$ and~$Ls$ are closed under~$\cdot_g$, that is, if and only
if $g(L\bar s)$ and~$g(Ls)$ are closed under the standard octonion
product.
\end{proposition}

\begin{proof}
Write $xy$ for $x \cdot_g y$, while $L\bar s$ and~$Ls$ remain the
lattices of~\eqref{eq:disp1}.  Since $g(L) = L$,
Lemma~\ref{lem:L-order} gives $L \cdot_g L = g^{-1}(L \cdot L)
\subseteq L$, so every block of a product of vectors of~$L^3$ lies
in~$L$.  For $u = (x, y, z) \in L^3$ put $a = x + y$, $b = x + z$, and
$w = x + y + z$.  As $y + z = a + b - 2x$ and $2L \subseteq L\bar s$
(Section~\ref{subsec:wilson-leech}), $u \in \LL$ exactly when $a, b
\in L\bar s$ and $w \in Ls$, and then $u = (a + b - w,\ w - b,\ w -
a)$.  For two such vectors $u$ and~$u'$, Definition~\ref{def:triple}
gives
\[
  P_1 + P_2 + P_3 = ww', \qquad
  P_1 + P_2 = ab' + ba' - aa' + 2(bb' - wb' + ww' - bw'),
\]
and the other two pairwise sums follow by the $\ZZ_3$ symmetry of the
routing.  If $L\bar s$ and~$Ls$ are closed, the block sum lies in~$Ls$
and, since $2L \subseteq L\bar s$, the pairwise sums lie in~$L\bar s$.
Conversely, $(a, 0, -a)$ and $(w, w, -w)$ lie in~$\LL$ for $a \in L\bar
s$ and $w \in Ls$; the product of two vectors of the first kind has $P_1
+ P_2 = -aa'$, and that of two of the second kind has block sum $ww'$,
so closure forces both inclusions.
\end{proof}

For $g = \sigma$, the lattices $\sigma(L\bar s) = 2\MO$ and $\sigma(Ls) =
p\MO$ of~\eqref{eq:multiples} are closed under the standard octonion
product, by the closure of~$\MO$ and Lemma~\ref{lem:sigmaLs-closed}, so
the proposition gives Theorem~\ref{thm:main} again.

\section{Related work}\label{sec:related}

\textbf{Wilson (2009)} \cite{Wilson2009} characterizes~$\LL$ as
$\{(x,y,z) \in L^3 : \text{conditions (1)--(3)}\}$, gives an explicit
listing of the 196{,}560 minimal vectors in the form
$3 \times 240 \times (1 + 16 + 16 \times 16)$, and proves that~$\LL$
is even and self-dual in his halved-norm convention.  The sublattice description is
the setting used in the present paper.

\textbf{Dixon (2010)} \cite{Dixon2010} obtains the integral octonions
from the octonion XY-product
\begin{equation}\label{eq:disp10}
  A \circ_{XY} B \;=\; (A X)(Y^{\dagger} B),
\end{equation}
where Dixon writes $X^{\dagger}$ for the conjugate~$\bar X$, and
$X, Y \in \OO$ are chosen so that $Y X^{\dagger}$ is a unit octonion
(the identity of the new multiplication).  When $Y^{\dagger} = X^{-1}$
the XY-product specializes to the X-product $A \circ_X B = (AX)(X^{-1}
B)$; Dixon uses this X-product with $X = \ell_0^{-1}$, where $\ell_0 =
\frac{1}{2}(e_0 + e_1 + \cdots + e_7)$, to obtain~$E_8$, and XY-products
for the laminated lattice~$\Lambda_{16}$.  For~$\LL$, Dixon writes the
minimal vectors of a copy of~$\LL$ (the componentwise conjugates of
Wilson's~\cite{Wilson2009}) in~$\OO^3$ with the standard octonion
product.\footnote{The construction
has an earlier precursor in Dixon~\cite{Dixon1995}, where a closely
related octonionic description of~$\LL$ is proposed.  A 2005 preprint
of Dixon's~\cite{Dixon2005} carried a formula whose third family of
minimal vectors required a small correction, which Wilson
\cite{Wilson2009} supplied in a bracketed note while independently
deriving the sublattice characterization.  The situation is similar
to the Kirmse--Bruck--Coxeter relationship recorded in
Appendix~\ref{app:history}: an early, structurally correct proposal
needing correction.}  We have independently reproduced this
construction.\footnote{\scriptref{leech\_dixon.py}
\scriptref{e8\_dixon.py}.}

\textbf{Baez--Egan (2014).}  In part~9 of the series \emph{Integral
octonions}~\cite{Baez2014}, Greg Egan constructs the Leech lattice from
the $E_8$ lattice alone; in the form of part~11, for sublattices
$L_1, L_2 \subset L$, each isometric to~$\sqrt{2}\,L$, with
$L_1 \cap L_2 = 2L$, the lattice
\begin{equation}\label{eq:egan}
  \LL(L_1, L_2) \;:=\; \bigl\{ (x,y,z) \in L^3 :
    x{+}y,\, x{+}z,\, y{+}z \in L_1;\; x{+}y{+}z \in L_2 \bigr\}
\end{equation}
is isometric to~$\LL$.  Wilson's lattice is $\LL = \LL(L\bar s, Ls)$
(Definition~\ref{def:wilson-leech}), and, with~$\sigma$ applied in each
block, $\sigma(\LL) = \LL(2\MO, p\MO)$ by~\eqref{eq:multiples} and
Lemma~\ref{lem:sigma-L}.  Applied in each block in the same way, an
isometry~$g$ of~$L$ with $g(L\bar s) = L_1$ and $g(Ls) = L_2$, which
exists by Appendix~\ref{app:mod2quotient}, carries $\LL$ onto
$\LL(L_1, L_2)$ and~$\star_g$ onto the untwisted product~$\star_{\mathrm{id}}$
of Section~\ref{sec:criterion}.  By Proposition~\ref{prop:criterion}, then,
$\star_{\mathrm{id}}$ closes on $\LL(L_1, L_2)$ exactly when $L_1$
and~$L_2$ are closed under the standard octonion product: it closes
on~$\sigma(\LL)$, which is Theorem~\ref{thm:main}, and not on~$\LL$
(Remark~\ref{rem:Ls-not-closed}).  There are $270$ sublattices of~$L$
isometric to~$\sqrt{2}\,L$, and each meets exactly $64$ of them in~$2L$,
so Egan's construction gives $17{,}280 = 270 \times 64$ lattices, all
distinct (part~11, and Egan's comments there~\cite{Egan2014}).
Appendix~\ref{app:mod2quotient} reads these pairs modulo~$2$, with the
product that $L/2L$ inherits from the standard octonion product.

In parts~9 and~10, Baez considers the 27-dimensional lattice
$\ZZ^3 \oplus L_L \subset \mathfrak{h}_3(\OO)$ of the exceptional Jordan
algebra, with integers on the diagonal and an off-diagonal Leech
lattice~$L_L$ at the standard scale, with entries in a maximal order of
integral octonions (here~$\MO$), so that $L_L \cong \LL/\sqrt{2}$.  Egan found~\cite{Egan2014} that this lattice
is \emph{not} closed under the Jordan product $x \circ y = xy + yx$ but
\emph{is} closed under the doubled product $2(xy + yx)$.  Egan further
enumerates at least $17{,}280$ distinct such lattices~$L_L$ for which
the associated 27-dimensional lattice is closed under this doubled
product: those of Egan's construction with~$\MO$ in place of~$L$ (see
part~10 for Egan's computational counts).  Baez leaves open the question
whether the 27-dimensional lattice of any such~$L_L$ is closed under the
undoubled Jordan product, and conjectures that the answer is negative.  We have likewise independently reproduced the full
Baez--Egan verification chain: Egan's simple roots of~$E_8$ and the two
sublattices built from them, the lemma that these meet in twice that
$E_8$ lattice, the Leech-lattice theorem for that pair, the
doubled-Jordan-product closure,\footnote{\scriptref{verify\_egan\_construction.py}.}
the derivation of the lower bound $244{,}035{,}421$ on the number of
sublattices of~$L^3$ isometric to~$\LL$,\footnote{An orbit-stabilizer bound, $\lceil 3!\,\lvert W(E_8)\rvert^3 / \lvert\Co_0\rvert\rceil$, with $W(E_8)$ the Weyl group of~$E_8$;
\scriptref{verify\_egan\_lower\_bound.py}.} and Egan's
count.\footnote{\scriptref{egan\_baez\_count.py}.}

\textbf{Anatomy of the Baez--Egan closure on $\LL$.}  The
Baez--Egan closure is on the 27-dimensional lattice
$\ZZ^3 \oplus L_L$, not on~$L_L$ alone.  At the paper's scale, for
off-diagonal Hermitian matrices $M_u, M_v \in \mathfrak{h}_3(\OO)$
representing $u = (x,y,z)$ and $v = (x',y',z')$ in $\LL \subset \OO^3$,
direct computation gives
\begin{equation}\label{eq:disp11}
  2(M_u M_v + M_v M_u) \;=\; \Delta \;+\; R,
\end{equation}
where $R$ is an off-diagonal residue lying in $(2L)^3 \subseteq \LL$
and $\Delta \in \ZZ^3$ is a diagonal piece whose entries are four times the
block-pair Euclidean inner products of $u$ and $v$:
\begin{equation}\label{eq:disp12}
  \Delta \;=\; 4\bigl(\langle x, x'\rangle + \langle y, y'\rangle,\;
                  \langle x, x'\rangle + \langle z, z'\rangle,\;
                  \langle y, y'\rangle + \langle z, z'\rangle\bigr),
\end{equation}
so $\operatorname{tr} \Delta = 8 \langle u, v\rangle$ recovers the
$\RR^{24}$ inner product of $u$ and $v$.  Projecting onto the
off-diagonal subspace defines a bilinear product $\varphi(u, v) := R$
under which $\LL$ is closed, so that $(\LL, +, \varphi)$ is an order
in the sense of Corollary~\ref{cor:order}.  For minimal vectors
$u \in \operatorname{Min}(\LL)$, squaring gives
$\operatorname{tr} \Delta = 8\Nrm(u) = 64$, so the projection is
nontrivial on every minimal-vector input.  The induced product
$\varphi$ is moreover commutative (the Jordan product is symmetric
in its arguments), whereas $\star$ is not
(Section~\ref{sec:identities}).  The products $\varphi$ and~$\star$
are therefore distinct bilinear products on~$\RR^{24}$ that both close
on~$\LL$.  Their spans differ as well.  Computed by the same
Hermite-normal-form method as for~$\star$
(Appendix~\ref{sec:span}),\footnote{\scriptref{verify\_phi\_span\_index.py}.}
the span index of~$\varphi$ is~$2^{25}$.  Unlike for~$\star$, the span
of~$\varphi$ does not contain~$2\LL$: its index exceeds
$[\LL : 2\LL] = 2^{24}$.

\textbf{Earlier precursors.}  Three works of 1982 and 1996
construct~$\LL$ from three copies of~$E_8$ or relate it to the
exceptional Jordan algebra.
John H.~Conway and Neil J.\,A.~Sloane~\cite{ConwaySloane1982} give 23
constructions of~$\LL$, one for each Niemeier lattice other than the
Leech lattice, including a construction from three copies of~$E_8$.  James Lepowsky
and Arne Meurman~\cite{LepowskyMeurman1982} describe~$\LL$ as a
sublattice of~$E_8^3$ containing~$2 E_8^3$ and generate the Conway
group~$\Co_0$ from a maximal $2$-local subgroup, which permutes the
three $E_8$ factors, and one extra automorphism.  Noam D.~Elkies and
Benedict H.~Gross~\cite{ElkiesGross1996} relate an integer form of the
exceptional Jordan algebra $\mathfrak{h}_3(\OO)$, with a positive
definite element of determinant~$1$ outside the orbit of the identity
matrix, to~$\LL$ and to the even lattice of rank~$26$ and
determinant~$3$ with no vectors of squared length~$2$.

\textbf{Leonard E.~Dickson (1923), Johannes Kirmse (1924), Kurt Mahler
(1942), Richard H.~Bruck and H.\,S.\,M.~Coxeter (1946).}  For the
history, see Appendix~\ref{app:history}.

\textbf{Daniele Corradetti (2026).}  Corradetti~\cite{Corradetti2026}
studies integral structures on~$E_8$ through the Okubo algebra,
introducing an integral order-3 automorphism of the octonions whose
Petersson isotope is the compact real Okubo algebra.  The
Coxeter--Dickson $E_8$ lattice then carries an octonionic, a
para-octonionic, and an Okubo integral structure, three multiplicative
systems that are not isomorphic.

\section{Conclusion}\label{sec:conclusion}

The Leech lattice~$\LL$ is an order under the $\ZZ_3$-symmetric
triple-octonion product of Definition~\ref{def:triple}
(Corollary~\ref{cor:order}), in Wilson's sublattice description.  The
proof (Section~\ref{sec:proof}) rests on one maximal order~$\MO$ of
integral octonions: the transposition~$\sigma$ fixes~$L$ and moves
Wilson's sublattices $L\bar s$ and~$Ls$ onto the left multiples
$2\MO$ and~$p\MO$, and the closure of~$\MO$ gives the inclusions the
proof needs.  More generally, for every twist~$g$, the product~$\star_g$
closes on~$\LL$ exactly when $L\bar s$ and~$Ls$ are closed
under~$\cdot_g$ (Proposition~\ref{prop:criterion}).  The order is
nonunital, has no nonzero idempotent, and is
not maximal; its automorphism group is $C_6 \times S_3$, and its span
index is~$2^{16}$, a value that depends on the twist
(Section~\ref{sec:properties} and Appendix~\ref{sec:span}).

Several questions remain open:
\begin{itemize}
\item Where $\Aut(\LL, +, \star)$ sits inside~$\Co_0$ up to conjugacy,
  and whether the order-6 automorphism of~$(\OO, \cdot_\sigma)$ generating the~$C_6$
  (Section~\ref{sec:autgroup}) admits an intrinsic description.
\item A classification of the algebra $(\RR^{24}, +, \star)$, which is
  nonassociative, nonalternative, nonflexible, non-power-associative,
  noncommutative, and nonunital (Section~\ref{sec:identities}).
\item A reason for the closure pattern of Appendix~\ref{sec:census}
  beyond the $2$-cycles, in particular for the failure of every
  $4$-cycle and for the exact one-half rates.
\end{itemize}

\section{Outlook}\label{sec:outlook}

A more natural classification of $(\RR^{24}, +, \star)$ may emerge
from a \emph{ternary} rather than binary viewpoint.  The Leech lattice
carries a prominent $\ZZ_3$ symmetry (the cyclic permutation of the
three octonion blocks), and~$\star$ is assembled from nine bilinear
building blocks routed by that $\ZZ_3$ action.  Composition algebras
with triality, as studied by
Elduque~\cite{Elduque1996,Elduque2000_Triality,Elduque2023_IsotropicNorm}
and others~\cite{SmithVojtechovsky2022,KamiyaOkubo2015,Okubo1995_Book},
are typically nonassociative, noncommutative, nonalternative, and
non-power-associative, the same negative classification $(\RR^{24},
+, \star)$ satisfies.  Yet they carry high internal symmetry, with
order-3 automorphisms playing the central role: the Okubo algebra has
automorphism group $\mathrm{SU}(3)/\ZZ_3$, and its product is built from an
order-3 automorphism of the octonions.  Jonathan I.~Hall's
monograph~\cite{Hall2019_Moufang} establishes that Moufang loops and
groups with triality are essentially the same object.  Recasting
$\star$ as a trace-like reduction of a ternary operation compatible
with the $\ZZ_3$-block symmetry is a natural direction for further
work.

Programs that model fundamental physics on normed-division-algebra,
triality, $S_3$ symmetry, or symmetric-composition-algebra structures
\cite{GresnigtGourlayVarma2023, Koeplinger2023,
FureyHughes2025_TrioOfTrialities, MarraniCorradettiZucconi2025}
are one setting in which a nonclassical bilinear product on~$\RR^{24}$
may be of interest.

\section*{Acknowledgments}

The author thanks Matthew Barley, Geoffrey M.~Dixon, Holger
P.~Petersson, Petr Vojt\v{e}chovsk\'{y}, and Robert A.~Wilson for
valuable discussion, advice, feedback, and corrections.

\section*{Declarations}

\textbf{Funding.}  The author received no funding for this work.

\textbf{Competing interests.}  The author has no competing interests
to declare.

\textbf{Data and code availability.}  All data and code are public:
the complete research record, including every verification script
cited in this paper, is available at the repositories of~\cite{repo}.

\textbf{Use of AI.}  AI assistance is a stated, material part of this
work, including the drafting of the manuscript text; the methodology,
the role of AI, and the human-designed verification and governance
process are documented in Appendix~\ref{sec:methodology}.  The author directed the work prompt
by prompt and reviewed the AI-assisted output; machine-obtained results
were verified under the protocol documented there.  The author takes
full responsibility for the content, including the accuracy of the
references.

\appendix
\section{The closure criterion modulo~\texorpdfstring{$2$}{2}}\label{sec:appendix-tables}\label{app:mod2quotient}

Because $L \cdot L \subseteq L$ (Lemma~\ref{lem:L-order}) and the
standard octonion product is $\ZZ$-bilinear, $(a + 2c) \cdot b \equiv a \cdot
b \pmod{2L}$ for $a, b, c \in L$, and likewise in the right factor.  The
product therefore descends to an $\FF_2$-bilinear product~$\bar\mu$ on
\begin{equation}\label{eq:disp14}
  \bar L := L/2L \cong \FF_2^8,
\end{equation}
the reduction of the order~$L$ modulo~$2$.  A lattice~$M$ with
$2L \subseteq M \subseteq L$ is the full preimage of~$M/2L$, so it is
closed under the standard octonion product exactly when $M/2L$ is a
subalgebra of~$(\bar L, \bar\mu)$.  For a twist~$g$
(Section~\ref{sec:criterion}), the lattices $g(L\bar s)$ and~$g(Ls)$
contain~$2L$, since $L\bar s \cap Ls = 2L$
(Section~\ref{subsec:wilson-leech}) and $g(2L) = 2L$.
Proposition~\ref{prop:criterion} therefore reads modulo~$2$: $\star_g$
closes on~$\LL$ if and only if $g(L\bar s)/2L$ and~$g(Ls)/2L$ are
subalgebras of~$(\bar L, \bar\mu)$.

\paragraph{Lagrangians and polarizations.}
The quadratic form
\begin{equation}\label{eq:disp17}
  q \colon \bar L \to \FF_2, \qquad q(v + 2L) \;:=\; \Nrm(v)/2 \pmod 2,
\end{equation}
is well defined because $L$ is even ($\Nrm(v + 2u) \equiv \Nrm(v)
\pmod 4$), and its associated bilinear form, the inner product
modulo~$2$, is nondegenerate because $L$ is unimodular.  A subspace
of~$\bar L$ is \emph{totally isotropic} when $q$ vanishes on it; it
then lies in its orthogonal complement, so its dimension is at
most~$4$.  The $4$-dimensional ones are the \emph{Lagrangians}, and a
\emph{polarization} of~$\bar L$ is an ordered pair of complementary
Lagrangians.

\paragraph{Wilson's pair and $\sigma$'s pair.}
The subspaces $L\bar s/2L$ and~$Ls/2L$ form a polarization.  They are
complementary, since $L\bar s + Ls = L$ and $L\bar s \cap Ls = 2L$;
each has dimension~$4$, since $[L : L\bar s] = [L : Ls] = 16$; and $q$
vanishes on both, since $\Nrm(\ell\bar s) = \Nrm(\ell s) =
2\Nrm(\ell) \in 4\ZZ$ for $\ell \in L$.  The isometry~$\sigma$ of~$L$
maps this polarization to
\begin{equation}\label{eq:disp15}
  V := \sigma(L\bar s)/2L = 2\MO/2L, \qquad
  W := \sigma(Ls)/2L = p\MO/2L,
\end{equation}
where the second equalities are~\eqref{eq:multiples}.\footnote{\scriptref{verify\_polarization\_lagrangians.py}.}
The identity twist does not close, since $Ls/2L$ is not a subalgebra
(Remark~\ref{rem:Ls-not-closed}).  The twist~$\sigma$ closes:
modulo~$2$, Lemma~\ref{lem:L-sigmaLsbar} says that $V$ is a left ideal
of~$(\bar L, \bar\mu)$, and Lemma~\ref{lem:sigmaLs-closed} that $W$ is
a subalgebra.\footnote{\scriptref{verify\_mod2\_quotient.py}.}
Moreover, $V$ is a two-sided ideal, since $2\MO$ is one in~$L$
(Remark~\ref{rem:Ls-not-closed}).\footnote{Hence
$(\bar L, \bar\mu)$ is not an octonion algebra in the
composition-algebra sense: those are simple.}

\paragraph{Which polarizations close.}
The Weyl group of~$E_8$, the isometry group of~$L$, acts on~$\bar L$
as the full orthogonal group of~$q$, with kernel~$\{\pm 1\}$
\cite[\S3.12.4]{WilsonBook2009}.  This group permutes the polarizations
transitively.  Indeed, each half of a polarization is its own orthogonal
complement, so the associated bilinear form pairs the two halves
$U_1$ and~$U_2$ nondegenerately: a basis $u_1, \ldots, u_4$ of~$U_1$
has a dual basis $u'_1, \ldots, u'_4$ in~$U_2$, and then
$q\bigl(\sum_i (a_i u_i + b_i u'_i)\bigr) = \sum_i a_i b_i$; the linear
map between two such bases preserves~$q$.  Hence every polarization is
$(g(L\bar s)/2L, g(Ls)/2L)$ for an isometry~$g$ of~$L$, and $\star_g$
closes on~$\LL$ exactly when both halves are subalgebras
of~$(\bar L, \bar\mu)$.

\paragraph{Relation to prior work.}
The orthogonal space~$(\bar L, q)$ is classical
\cite[\S3.12.4]{WilsonBook2009}.  Every sublattice~$M$ of~$L$
isometric to~$\sqrt2\,L$ contains~$2L$, since its dual lattice,
$\tfrac12 M$, contains~$L$.  Lepowsky and Meurman
\cite[\S2]{LepowskyMeurman1982} identify the lattices between~$2L$
and~$L$ that are isometric to~$\sqrt2\,L$ with the preimages of the
Lagrangians, and construct a lattice isometric to~$\LL$ from any
polarization; in their introduction they note that work of Jacques Tits
contains what amounts to the same construction.  Egan's construction (Section~\ref{sec:related}) is
of this form.  In the Baez series, Tim Silverman arrived at the
identification of these sublattices with the Lagrangians independently
(part~12 of~\cite{Baez2014}), read Egan's construction modulo~$2$, and
found that Wilson's pair is a polarization~\cite{Silverman2016}.  David
Calderbank related Egan's construction to that of Lepowsky and
Meurman~\cite{Calderbank2023}.
The present appendix adds the product~$\bar\mu$ inherited from the
standard octonion product, under which $V$ and~$W$ are a two-sided
ideal and a subalgebra.

\section{Historical note: integral Cayley numbers,
1923--1946}\label{app:history}

The integral octonions this paper builds on were developed across four
primary papers between 1923 and 1946.  Each was obtained in the
original, and its key computational content re-derived in exact
arithmetic in this project (for Mahler's covering-radius theorem, by
exact-arithmetic sampling).\footnote{\scriptref{verify\_dickson\_1923.py}
\scriptref{verify\_kirmse\_1924.py}
\scriptref{verify\_kirmse\_1924\_forensic.py} \scriptref{verify\_mahler\_1942.py}
\scriptref{verify\_coxeter\_1946.py}.  Verification notes accompany each,
in \texttt{paper/reviews/}.}

\noindent\textbf{Dickson (1923)} \cite{Dickson1923} gave the first
construction: an explicit 8-element basis of an order~$O_D$ for the
Cayley--Dickson doubled product, with all 64 basis products
closing.  His Theorem~XV states three maximal orders containing the
eight basis units; the true count is \emph{seven}.  The cause is in his derivation (p.~321): he assumes
the orders contain the Hurwitz quaternions, and exactly three of the
seven do.

\noindent\textbf{Kirmse (1924)} \cite{Kirmse1924} took up the problem
independently.  His multiplication table is itself a correct
octonion algebra, and he derived a correct \emph{necessary} form
for the maximal orders, without proving sufficiency; thirty
$E_8$-type lattices realize the form.  Among the eight of those lattices that
are stable under multiplication by the basis units (an unstated
narrowing of the criterion, in our reconstruction), he exhibited the
one, $J_1$ (p.~70), that is not an order.\footnote{Kirmse's 1928
corrigendum~\cite{Kirmse1928} corrects an interchange of his left and
right multiplication matrices and his definition of ideals; it does not
concern the maximal orders.}  On this reading, the
remaining seven of his eight are the genuine maximal orders.  Kirmse also develops an
ideal-theoretic treatment of the integral octonions, a substantive
companion contribution that is taken up again in Mahler (1942),
and methodically takes algebra alternativity as a building
block, anticipating the structural role alternativity plays in
later integral-octonion arithmetic.

\noindent\textbf{Mahler (1942)} \cite{Mahler1942} developed the
arithmetic of Dickson's order~$O_D$: for every Cayley number~$X$ there is
an integral~$G$ with $\Nrm(X-G)\le\tfrac{15}{16}$ (sharper than
Dickson's $\tfrac54$; Mahler states without proof that $\tfrac12$ is the
best constant), and every left or right ideal is principal, generated by
an integer multiple of an integral Cayley number of norm $1$, $2$, or $4$.

\noindent\textbf{Coxeter (1946)} \cite{Coxeter1946} identifies
Kirmse's error with an explicit nonclosure witness and records that
the remedy is due to Bruck: ``Bruck's domain $J$ can be
derived from $J_1$ by transposing two of the $i$'s'' (\S4).  This is
a transposition of two imaginary basis units, a linear involution
of~$\RR^8$; the $\binom72 = 21$ such transpositions yield the seven
maximal orders in seven sets of three.  Coxeter's \S\S6--13 develop the
geometry of the corrected order as the $E_8$ lattice, the vertex set of
the honeycomb~$5_{21}$ in his notation.  A postscript (\S14)
records his after-the-fact awareness of Dickson (1923) and Mahler (1942) and
identifies Dickson's undercount with exactly the cause above.

It is curious to find that the Coxeter--Bruck transposition would
prove useful here: the same linear involution of~$\RR^8$, acting
now on a different object.  The Coxeter--Bruck transposition acts
on a maximal-order candidate that sits asymmetrically in~$\RR^8$;
the present~$\sigma$ fixes the coordinate-symmetric placement~$L$
and does its work on the sublattice~$Ls$ of Wilson's Leech
construction.

\section{Research methodology}\label{sec:methodology}

This research was conducted as a collaboration between a human
researcher and an AI agent (Claude, Anthropic) over a series of
directed prompts, governed by a written methodology (the
``Manifesto'') established at the outset.

The human researcher (the author) brought the initial hypothesis
(motivated by prior work on Okubo algebras and autotopies of normed
composition algebras~\cite{Koeplinger2023}), established the research methodology, and
maintained control of the inquiry through prompt-by-prompt review.
The AI agent performed the computational work (implementing octonion
algebras, Wilson's Leech lattice construction, trial algebras, and
closure tests) and contributed mathematical analysis.  It also drafted
the manuscript text and made every revision to it, at the author's
prompt-by-prompt direction.

Structural proposals from a second AI system (the Hermes agent, Nous
Research, running DeepSeek V4 Flash through OpenCode) prompted the
Lagrangian framing of Appendix~\ref{app:mod2quotient}; the AI agent
(at that stage Claude Opus~4.7) verified and corrected them.\footnote{The proposals
verbatim, the verification, and the corrections: \cite{repo}~%
\href{https://github.com/koeplinger/leech\_alg/blob/main/paper/hermes\_mod2\_structure\_supplement.pdf}%
{\nolinkurl{paper/hermes_mod2_structure_supplement.pdf}};
\scriptref{verify\_discovery1\_W\_subalgebra.py}
\scriptref{verify\_discovery2\_V\_isotropic.py}.}

Tests were written from first principles to independently verify
every foundation.  This independent verification is how we know that the
properties stated in the literature are correctly implemented.
Closure of the twisted product was first observed on over
12{,}000{,}000 random pairs of minimal vectors with zero failures,
which prompted the search for a symbolic proof.

\paragraph{Symbolic proof and numeric validation.}
The proof of Theorem~\ref{thm:main} given in Section~\ref{sec:proof}
is symbolic.  Its finite content, the closure of~$\MO$ on a
$\ZZ$-basis and the identities~\eqref{eq:multiples} on eight products
each, is verified in exact rational
arithmetic.\footnote{\scriptref{verify\_structural\_proof.py}.}
No floating-point or sampled arithmetic enters the proof.
Independently, the four lemmas themselves were verified on
$\ZZ$-bases: Lemma~\ref{lem:sigma-L} on the eight basis vectors of~$L$,
and Lemmas~\ref{lem:L-order} to~\ref{lem:sigmaLs-closed} basis pair by
basis pair.  This was done in two separate computer-algebra
implementations, a Python implementation using
\texttt{fractions.Fraction} exact rationals and a GAP/LOOPS
implementation~\cite{GAP, LOOPS} using GAP's native exact rational
arithmetic.\footnote{\scriptref{symbolic\_proof\_checks.py}
\scriptref{verify\_structural\_proof.py} \gapref{tests/test\_lemmas.g}.}  This redundancy is validation, not
proof: it serves to detect implementation errors, but the logical
content of the proof rests on Section~\ref{sec:proof}.

The methodology required:
\begin{itemize}
\item Every claim backed by a reference or computation.
\item Every prompt logged, with actions recorded.
\item Every conclusion independently verified.
\item Forward-evolving corrections: the prompt logs are never rewritten,
  and later entries correct earlier errors; other documents are corrected
  in place, date stamped with the update, with their history in git.
\end{itemize}

The path to the result was not direct.  Six trial algebra families
were explored and ruled out before the transposition twist was
identified in trial~007.  The systematic methodology was essential:
understanding the anatomy of failure (which conditions fail, for
which vector types, and why) before trying the next approach.  The failure analysis of the untwisted product (Wilson's
condition~(3) per Definition~\ref{def:wilson-leech}, exclusively
Type~3$\times$Type~3) directly pointed to the
mechanism that the twist needed to correct.

Iterative convergence toward a robust result, including learning
along the way what it is we are actually doing, is how study,
research, and discovery have always proceeded.  Mistakes were made
on both sides: the AI agent occasionally departed from the agreed
protocol or produced incorrect intermediate analyses, and the human
author at times framed prompts under misconceptions or accepted
intermediate results without sufficient scrutiny.  In each case, the
error was caught by the systematic cross-checking that the
methodology requires (explicit verification, prompt-by-prompt
review, and independent reimplementation in Python and GAP) and
corrected forward.  What is novel here is the
ability to trace the entire path, beginning to end, through the
complete list of prompts.

The full record is available at the repositories of~\cite{repo}: all
prompts, code, evidence, and git history.  The reader can form
their own assessment of what this means for questions of discovery
and attribution.

\section{Identities and product norms on minimal
vectors}\label{app:samples}

This appendix estimates, by sampling, how often the identities of
Section~\ref{sec:identities} from left alternativity to symmetric
composition hold on minimal vectors, and records the norms of products
of minimal vectors.  Samples were drawn uniformly from
$\operatorname{Min}(\LL) \times \operatorname{Min}(\LL)$ for the
binary identities, from $\operatorname{Min}(\LL)$ for the cube and
quartic identities, and from $\operatorname{Min}(\LL)^3$ for the
symmetric-composition identity.  Sample size is $n = 10^6$ per
identity; exact integer arithmetic
throughout.\footnote{\scriptref{verify\_section5\_properties.py}.}  The
standard error is $\mathrm{SE} = \sqrt{f(1-f)/n}$, where $f$ is the
fraction of samples on which the identity holds.

\begin{center}
\begin{tabular}{@{}ll@{}}
\toprule
Identity & Fraction of samples (\% $\pm$ SE) \\
\midrule
Left alternativity & $(0.188 \pm 0.0043)\%$ \\
Right alternativity & $(0.188 \pm 0.0043)\%$ \\
Flexibility & $(8.37 \pm 0.028)\%$ \\
Cube power-associativity & $(30.70 \pm 0.046)\%$ \\
Quartic power-associativity & $(12.75 \pm 0.033)\%$ \\
Commutativity & $(0.018 \pm 0.0014)\%$ \\
Norm multiplicativity & $(46.77 \pm 0.050)\%$ \\
Symmetric composition & $(25.37 \pm 0.044)\%$ \\
\bottomrule
\end{tabular}
\end{center}

On the $10^6$ samples drawn for the cube and quartic identities, every
sample that satisfied the quartic identity also satisfied the cube
identity.

The self-product norms $\Nrm(u \star u)$ of all 196{,}560 minimal
vectors take exactly the values $\{16, 32, 48, 64, 80, 96,
128\}$.\footnote{\scriptref{verify\_idempotents\_min\_shell.py}.}  On
$10^6$ sampled pairs, the pair-product norms $\Nrm(u \star v)$ took
only the values $\{16, 32, 48, 64, 80, 96, 112, 128\}$, all multiples
of~16, with mode~64 (the norm-multiplicative value $\Nrm(u)\Nrm(v) = 8
\times 8$).\footnote{\scriptref{verify\_pair\_norm\_histogram.py}.}

\section{Other twists: closure and span of the
image}\label{sec:census}\label{sec:span}

This appendix replaces~$\sigma$ by other twists
(Section~\ref{sec:criterion}) and tests two properties of the
resulting triple product: whether it still closes on~$\LL$, and how
much of~$\LL$ its products reach.

When $\star_g$ closes on~$\LL$, the integer combinations of its
products $u \star_g v$ ($u, v \in \LL$) form a sublattice of~$\LL$,
the \emph{span} of~$\star_g$.  The \emph{span index} of~$\star_g$ is
the number of disjoint translates of the span that make up~$\LL$.  Span index~$1$ means that every vector
of~$\LL$ is an integer combination of products.  An order that
contains a multiplicative identity always has span index~$1$, since
every vector $u = u \cdot 1$ is a product; the ambient algebra here
has no multiplicative identity, and the span index can be larger.

Both properties are decided by finitely many products.  By
bilinearity it suffices to take the $24 \times 24 = 576$ products of
pairs of vectors from a $\ZZ$-basis of~$\LL$.  The product~$\star_g$
closes on~$\LL$ if and only if all of them lie in~$\LL$, and the span
index is the index of the lattice they span, computed exactly by
Hermite normal form.\footnote{\scriptref{verify\_product\_span\_index.py}.}

As an experiment, we generalize from the transposition~$\sigma$ to all
permutations~$\pi$ of the seven imaginary basis units
$e_1, \ldots, e_7$, each fixing~$e_0$.  These form the symmetric
group~$S_7$, with~$\sigma$ as a special case.  For every nonidentity cycle type
of~$S_7$ we find:\footnote{\scriptref{verify\_all\_permutations\_exact.py}
\scriptref{verify\_twist\_span\_index\_census.py}.}
\begin{center}
\begin{tabular}{@{}lrrrr@{}}
\toprule
Cycle type of $\pi$ & Permutations & Close & Rate & Span index \\
\midrule
$2$-cycle     & $21$  & $21$  & $100\%$  & $2^{16}$ \\
$3$-cycle     & $70$  & $35$  & $50\%$   & $2^{10}$ \\
$(2,2)$       & $105$ & $42$  & $40\%$   & $2^{10}$ \\
$4$-cycle     & $210$ & $0$   & $0\%$    & n/a \\
$(3,2)$       & $420$ & $126$ & $30\%$   & $2^{16}$ \\
$5$-cycle     & $504$ & $252$ & $50\%$   & $2^{10}$ \\
$(2,2,2)$     & $105$ & $21$  & $20\%$   & $2^{16}$ \\
$(3,3)$       & $280$ & $112$ & $40\%$   & $2^{10}$ \\
$(4,2)$       & $630$ & $294$ & $46.7\%$ & $2^{10}$ \\
$6$-cycle     & $840$ & $210$ & $25\%$   & $2^{16}$ \\
$(3,2,2)$     & $210$ & $105$ & $50\%$   & $2^{10}$ \\
$(4,3)$       & $420$ & $168$ & $40\%$   & $2^{16}$ \\
$(5,2)$       & $504$ & $42$  & $8.3\%$  & $2^{16}$ \\
$7$-cycle     & $720$ & $336$ & $46.7\%$ & $2^{10}$ \\
\midrule
Total         & $5{,}039$ & $1{,}764$ & $35.0\%$ & \\
\bottomrule
\end{tabular}
\end{center}
The $2$-cycles, the transpositions of Definition~\ref{def:twist},
close by Theorem~\ref{thm:main}.  By Proposition~\ref{prop:criterion},
the identity permutation does not close, since $Ls$ is not closed
under the standard octonion product (Remark~\ref{rem:Ls-not-closed}).
Each permutation~$\pi$ preserves~$L$, as $\sigma$ does
(Lemma~\ref{lem:sigma-L}), so by the same proposition the column
``Close'' counts the permutations~$\pi$ for which $\pi(L\bar s)$
and~$\pi(Ls)$ are closed under the standard octonion product.  For the
$3$-cycles, exactly one of the two cyclic orientations of each
three-element subset of imaginary indices
closes.\footnote{\scriptref{verify\_twist\_characterization.py}.}  No
comparable pattern is given here for the other types that close in part,
nor a reason for the failure of every $4$-cycle.

The span index of a closing permutation depends only on its parity:
it is $2^{16}$ for every odd one and $2^{10}$ for every even one (last
column).

\begin{proposition}\label{prop:span}
The span of~$\star$ consists of the vectors $(w, w, w) + 2u$ with $w \in
Ls$ and $u \in \LL$.  In particular, $\star$ has span index $2^{16} =
65{,}536$, and its span contains~$2\LL$.\footnote{\scriptref{verify\_product\_span\_index.py}
\scriptref{verify\_product\_span\_structure.py}
\scriptref{verify\_product\_span\_proof.py}.}
\end{proposition}

\begin{proof}
Write~$S$ for the span.  By~\eqref{eq:multiples} and the closure
of~$\MO$, $L\bar s \cdot_\sigma L\bar s \subseteq 2L\bar s$ and $L\bar s
\cdot_\sigma Ls + Ls \cdot_\sigma L\bar s \subseteq L\bar s$.  The
products $(yp)z$ with $y, z \in \MO$ span~$\MO$, and so do the products
$y(pz)$.  Indeed, $v := -\tfrac12(e_1 + e_2)s = \tfrac12(e_0 + e_2 + e_3 -
e_7)$ lies in~$\MO$ by~\eqref{eq:MO}, as $s \in L$, and $pv + vp = 2v_0\,p
- 2\langle p, v\rangle e_0 = p - e_0$, as $p$ is purely imaginary.  The
right alternative law $(yx)x = yx^2$ with $x = p + v$ gives $(yp)v + (yv)p
= y(pv + vp) = yp - y$, so $y = (yp)e_0 - (yp)v - \bigl((yv)p\bigr)e_0$
for every $y \in \MO$; the left alternative law gives $y = e_0(py) - v(py)
- e_0\bigl(p(vy)\bigr)$ in the same way.  Now $\sigma$ maps
$Ls \cdot_\sigma Ls$ onto the products $(py)(pz) = p\bigl((yp)z\bigr)$
of~\eqref{eq:moufang-imaginary}, and $L\bar s \cdot_\sigma Ls$ onto the
products $2y(pz)$, with $y, z \in \MO$.  Hence the products in $Ls
\cdot_\sigma Ls$ span~$Ls$, and those in $L\bar s \cdot_\sigma Ls$
span~$L\bar s$.

From here on, write $xy$ for $x \cdot_\sigma y$, while $L\bar s$ and~$Ls$
remain the lattices of~\eqref{eq:disp1}.  By the proof of
Proposition~\ref{prop:criterion}, $\LL$ is spanned by the vectors $(a, 0,
-a)$, $(b, -b, 0)$ with $a, b \in L\bar s$ and $(-w, w, w)$ with $w \in
Ls$ (their sum is the vector $(a + b - w, w - b, w - a)$ there).
Definition~\ref{def:triple} gives
\begin{align*}
  (a, 0, -a) \star (a', 0, -a') &= (aa', -2aa', aa'), &
  (b, -b, 0) \star (b', -b', 0) &= (bb', bb', -2bb'), \\
  (a, 0, -a) \star (b, -b, 0) &= (2ab, -ab, -ab), &
  (a, 0, -a) \star (-w, w, w) &= (-2aw, 2aw, 0), \\
  (b, -b, 0) \star (-w, w, w) &= (-2bw, 0, 2bw), &
  (-w, w, w) \star (-w', w', w') &= (3ww', -ww', -ww'),
\end{align*}
and the products in the opposite order are the same with the octonion
factors on the right exchanged.  For $c \in L\bar s$ and integers $k_1,
k_2, k_3$ with $k_1 + k_2 + k_3 = 0$, the vector $(k_1 c, k_2 c, k_3 c)$
lies in~$\LL$ (Definition~\ref{def:wilson-leech}).  So, by the inclusions
above, every product other than the products of two vectors $(-w, w, w)$
lies in~$2\LL$.  These are $2(2c, 0, 0) - (c, c, c)$ with $c = ww' \in
Ls$, and $(2c, 0, 0) \in \LL$ as $2L \subseteq L\bar s \cap Ls$.  Hence,
by bilinearity, every vector of~$S$ has the stated form.

Conversely, $S$ contains $(3c, -c, -c)$ for every $c \in Ls$, since the
products in $Ls \cdot_\sigma Ls$ span~$Ls$.  The cyclic shift of the
blocks preserves~$\LL$ and, by the $\ZZ_3$ symmetry of the routing,
commutes with~$\star$, so it maps~$S$ onto itself, and $S$ also contains
$(-c, 3c, -c)$ and $(-c, -c, 3c)$.  Adding the three gives $(c, c, c) \in
S$, hence $(4c, 0, 0) \in S$ and $2(-c, c, c) = 2(c, c, c) - (4c, 0, 0)
\in S$.  Since the products in $L\bar s \cdot_\sigma Ls$ span~$L\bar s$,
the products of $(a, 0, -a)$ and of $(b, -b, 0)$ with $(-w, w, w)$ give
$2(c, -c, 0)$ and $2(c, 0, -c)$ in~$S$ for every $c \in L\bar s$.  So $S$
contains twice the vectors $(a, 0, -a)$, $(b, -b, 0)$ and $(-w, w, w)$
that span~$\LL$, hence~$2\LL$, and with the vectors $(c, c, c)$, $c \in
Ls$, every vector of the stated form.

Finally, if $(c, c, c) \in 2\LL$ with $c \in Ls$, then $\tfrac32 c \in
Ls$ by condition~(3), so $\tfrac12 c \in Ls$.  Hence $S/2\LL \cong
Ls/2Ls$ has order~$2^8$, and $[\LL : S] = 2^{24}/2^8 = 2^{16}$.
\end{proof}

The transposition~$\sigma$ is also the reflection in the root
$e_1 - e_2$ of~$L$, and reflections in the other roots of~$L$ are
twists as well.  For the reflection~$\rho$ in the root
$r = \tfrac12(e_0 + e_1 + \cdots + e_6 - e_7)$, the triple
product~$\star_\rho$ closes on~$\LL$ with span
index~$1$.\footnote{\scriptref{verify\_twist\_span\_index\_census.py}.}

The span index is invariant under ring isomorphism, since an
isomorphism maps products to products and hence the span onto the
span.  Orders with different span indices are therefore not
isomorphic.  By the parity split, an order $(\LL, +, \star_\pi)$ with
$\pi$ odd (span index $2^{16}$, as for~$\sigma$) is not isomorphic
to one with $\pi$ even (span index $2^{10}$), and neither is
isomorphic to $(\LL, +, \star_\rho)$ (span index~$1$).


\end{document}